\documentclass[12pt,twoside]{article}
\usepackage{amsfonts}
\usepackage{amssymb,amsmath}
\usepackage{float}
\newtheorem{theorem}{Theorem}

\newfloat{Exhibit}{}{}
\floatname{Exhibit}{Exhibit}

\begin{document}
\begin{center}
\textbf{Mathematical Structures defined by
Identities III}
\bigskip
\\Constantin M. Petridi
\\ cpetridi@math.uoa.gr
\\ cpetridi@hotmail.com
\end{center}
\par
\vspace{2pt}
\begin{center}
\small{
\begin{tabular}{p{11cm}}
\textbf{Abstract.}We extend the theory (\textbf{formal part
only}) of algebras with one binary operation (our paper
arXiv:math/0110333v1 [math.RA] 31 Oct 2001) to algebras with
several operations of any arity.
\bigskip
\end{tabular}
}
\end{center}

\vspace{10pt} \textbf{1 Introduction}\par\vspace{10pt}
We refer to our papers \cite{Petridi1} and \cite{Petridi2} for concepts, notations, definitions, notes and remarks.
\par\vspace{10pt}In subsection $4.1$ of \cite{Petridi1} we briefly outlined our ideas of generalizing the method of tableaux to algebras with several operations
\[
V_{1}(x_1,\,x_2,\,\cdots,\,x_{\alpha_1}),\,V_{2}(x_1,\,x_2,\,\cdots,\,x_{\alpha_2}),\,\cdots,\,V_{k}(x_1,\,x_2,\,\cdots,\,x_{\alpha_k})
\]
satisfying axiomatically defined identities and indicated the way
of how to proceed. The project is now carried out. The technique
applied is the same as in Formal Part of \cite{Petridi1}. The
crucial fact that the number $I_{n}^{V_{1}V_{2}\cdots V_{k}}$ of
\textbf{formally} reducible identities can be calculated by
exactly the same method used for $I_{n}^{V_1}(=I_{n})$ seems to
hold true. Algebras with only binary operations are discussed.
For algebras with two binary operations $V(x,\,y)$ and $W(x,\,y)$
the proof is given in detail.
\par\vspace{10pt}Algebras with operations of any arity can be treated by reduction to a well defined set of algebras with binary operations.
\par\vspace{10pt}Research and exposition of the general theory are impeded by problems of construction and inspection of the tableaux $T_{n}$ whenever $n$ is greater than $3$. This is due to the fast growth of the Catalan numbers $(S_{n}\sim\dfrac{4^n}{\pi^{\frac{1}{2}}n^{\frac{3}{2}}})$ and their generalizations, let alone problems of printing and publication. Programs designed to seek the structures resulting from a given identity failed after a few steps (blow-ups). Exposition therefore is limited to illustrate the theory on the worked example of tableau $T_{3}$.
\par\vspace{10pt}Still, the concrete new findings reached in this case corroborate further our fundamental thesis that there is a scarcity of existing mathematical structures in the sense that the frequency of \textbf{irreducible} identities goes to zero with increasing $n$. Seen historically, this also explains why mathematics, in the course of time, has developed the way it did with associativity $V(V(x,\,y),\,z)=V(x,\,V(y,\,z))$, the simplest structure, reigning supreme over the mathematical landscape.
All other essential mathematical structures, found or created by
research such as e.g. Groups, Fields, Vector Spaces, Lie
Algebras, etc, ... include in their axiom system (signature) at
least one binary operation obeying the law of associativity.
\par\vspace{10pt} We conclude with a note on the connection with
Formal Languages.
\par\vspace{20pt}

\textbf{2 Operations and their Iterates}\par\vspace{10pt} Given
$k$ operations
$V_{1}^{\alpha_{1}}(x_1,\,x_2,\,\cdots,\,x_{{\alpha}_{1}}) $ of
arity $\alpha_{1}$,
$V_{2}^{\alpha_{2}}(x_1,\,x_2,\,\cdots,\,x_{{\alpha}_{2}}) $ of
arity $\alpha_{2}$,$\ldots$,
$V_{k}^{\alpha_{k}}(x_1,\,x_2,\,\cdots,\,x_{{\alpha}_{k}}) $ of
arity $\alpha_{k}$, their $n-$iterates containing the operation
$V_{1}^{{\alpha}_{1}}$ $p_1-$times, the operation
$V_{2}^{{\alpha}_{2}}$ $p_2-$times,$\cdots$, the operation
$V_{k}^{{\alpha}_{k}}$ $p_k-$times are symbolized by \vspace{10pt}
\[ J_{i}^{n}\begin{pmatrix} V_{1}^{{\alpha}_{1}} & V_{2}^{{\alpha}_{2}} & \cdots & V_{k}^{{\alpha_{k}}} \\ p_1 & p_2 & \cdots & p_k \end{pmatrix},\,\,\alpha_{i}\geq{0},\,\,p_i\geq{0}.\]
\par\vspace{10pt}
 The order of the iterate is

\[ n=p_1+p_2+\cdots+p_k \]
\par\vspace{10pt}
and the number of its variable places is

\[ (\alpha_1-1)p_1+(\alpha_2-1)p_2+\cdots+\alpha_k(p_k-1)+1.\]

\par\vspace{10pt}The index $i$ runs from $1$ to $S_{n}^{V_1\,\cdots\,V_k}$.
We call $S_{n}^{V_1\,\cdots\,V_k}$ the \textbf{n-th Catalan
number} of the structure.
\par\vspace{10pt}The numbers $S_{n}^{V_1\,\cdots\,V_k}$ are the Taylor coefficients, at $t=0$, of the formal generating function
\[ \phi_{V_1\,\cdots\,V_k}(t)=\sum_{n=0}^{\infty}S_{n}^{V_1\,\cdots\,V_k}\,\,t^n. \]

The sequence $S_{n}^{V_1\,\cdots\,V_k}$ can be calculated recursively from

\begin{align*} S_{n+1}^{V_1\,\cdots\,V_k}=\sum_{\begin{subarray}{c} x_1+x_2+\cdots+x_{\alpha_1}=n \\
                                                         x_{i}\geq0\end{subarray}}
{S_{x_{1}}^{V_1\,\cdots\,V_k}\,S_{x_{2}}^{V_1\,\cdots\,V_k}\cdots S_{x_{\alpha_1}}^{V_1\,\cdots\,V_k}}\,+\\
\sum_{\begin{subarray}{c} x_1+x_2+\cdots+x_{\alpha_2}=n \\
                                                         x_{i}\geq0\end{subarray}}
{S_{x_{1}}^{V_1\,\cdots\,V_k}\,S_{x_{2}}^{V_1\,\cdots\,V_k}\cdots S_{x_{\alpha_2}}^{V_1\,\cdots\,V_k}}\,+ \\
\cdots\cdots\cdots\cdots\cdots\cdots\cdots\cdots\cdots\cdots \cdots\cdots \\
\sum_{\begin{subarray}{c} x_1+x_2+\cdots+x_{\alpha_k}=n \\
                                                         x_{i}\geq0\end{subarray}}
{S_{x_{1}}^{V_1\,\cdots\,V_k}\,S_{x_{2}}^{V_1\,\cdots\,V_k}\cdots S_{x_{\alpha_k}}^{V_1\,\cdots\,V_k}}.
\end{align*}

According to E. Catalan\footnote{See L.E. Dickson: History of the
theory of numbers, Vol.2} the number of solutions of the
Diophantine equation $x_1+x_2+\cdots+x_{\alpha_i}=n$ is
$\binom{a_{i}+n-1}{n}.$
\par\vspace{10pt}The function $\phi_{V_1\,\cdots\,V_k}(t)$ is a solution of the functional equation

\[ \dfrac{\phi_{V_1\,\cdots\,V_k}(t)-1}{t}=\sum_{i=1}^{k}\,\,(\phi_{V_1\,\cdots\,V_k}(t))^{\alpha_{i}}. \]

with the initial condition $\varphi(0)=1$.
\par\vspace{10pt}For $k=1,\,\alpha_1=2$ we get the classical Catalan numbers $S_{n}=\frac{1}{n+1}\binom{2n}{n}$, which count the $n-$iterates (parenthesizing) of $V_{1}(x,\,y)$. Their recursion formula is

\[ S_{n+1}^{V_{1}}=\sum_{\begin{subarray}{c} x_1+x_2=n \\ x_{i}\geq0 \end{subarray}}S_{x_1}^{V_1}\,S_{x_2}^{V_1} \]

and the functional equation becomes

\[ \dfrac{\phi_{V_1}(t)-1}{t}=(\phi_{V_1}(t))^{2}, \]

giving

\[
\phi_{V_1}(t)=\dfrac{1-\sqrt{1-4t}}{2t}=\sum_{n=0}^{\infty}\,\,S_nt^n.
\]

\par\vspace{10pt}For $k=1,\,\alpha_1=\alpha$ we obtain the higher Catalan numbers $\frac{1}{(\alpha-1)n+1}\binom{\alpha n}{n} $ whose generating function $\phi_{\alpha}(t)$ satisfies
\[
\dfrac{\phi_{\alpha}(t)-1}{t}=(\phi_{\alpha}(t))^{\alpha}.
\]
\par\vspace{20pt}

\textbf{3 Binary Operations}\par\vspace{10pt}
We will now examine the case of two binary operations $V(x,\,y)$ and $W(x,\,y)$. Since $k=2,\,\alpha_1=\alpha_2=2$ the corresponding generating function which gives the number of iterates of order $n$ is

\[ \dfrac{\phi_{VW}(t)-1}{t}=2(\phi_{VW}(t))^{2},\,\,\phi(0)=1. \]

Solving the quadratic  equation we obtain

\[ \phi_{VW}(t)=\dfrac{1-\sqrt{1-8t}}{4t}=\sum_{n=0}^{\infty}\,2^n\,S_n\,t^n, \]

where $S_{n}=\frac{1}{n+1}\binom{2n}{n}$ are the ordinary Catalan
numbers.
\par\vspace{10pt}Hence the number of iterates of order $n$ is $S_{n}^{VW}=2^nS_n,\;\;n\geq 1$,
the first of which are
\[ \begin{array}{c} n \\ \hline \\ 1 \\ 2 \\ 3 \\ 4 \\ 5 \\ \vdots
 \end{array}\,\,\,\begin{array}{c} S_{n}^{VW} \\ \hline \\ 2 \\ 8 \\ 40 \\ 224 \\ 1344 \\ \vdots \end{array}.\]

Following the same rules of formation as done in \cite{Petridi1}, the first three $A-$tableaux are

\[ \begin{array}{c}
\underline{\text{$T_{1}$}}\\\\
Vxx\\
Wxx \\ \\ \\
\end{array}\,\,\,\,\,\,\,\,\,\,\,\,\,\,\,\,\,
\begin{array}{c}
\underline{\text{$T_{2}$}}\\\\
VVxxx\,\,\,WVxxx\\
VWxx\,\,\,WWxxx\\
VxVxx\,\,\,WxVxx\\
VxWxx\,\,\,WxWxx
\end{array} \]

\begin{center} $\underline{T_{3}}$ \end{center}\vspace{-0.9cm}
\begin{minipage}[h]{2.5 cm}
\[ \begin{array}{c}
VVVxxxx\\
VVWxxxx\\
VVxVxxx\\
VVxWxxx\\
VVxxVxx\\
VVxxWxx
\end{array} \]
\end{minipage}
\begin{minipage}[h]{2.5 cm}
\[ \begin{array}{c}
WVVxxxx\\
WVWxxxx\\
WVxVxxx\\
WVxWxxx\\
WVxxVxx\\
WVxxWxx
\end{array} \]
\end{minipage}
\begin{minipage}[h]{2.5 cm}
\[ \begin{array}{c}  \\ \cdots \\ \cdots \\ \cdots \\ \cdots \\ \cdots \\ \cdots\end{array} \]
\end{minipage}
\begin{minipage}[h]{2.5 cm}
\[ \begin{array}{c}
VVxxWxx\\
VWxxWxx\\
VxWVxxx\\
VxWWxxx\\
VxWxVxx\\
VxWxWxx
\end{array} \]
\end{minipage}
\begin{minipage}[h]{2.5 cm}
\[ \begin{array}{c}
WVxxWxx\\
WWxxWxx\\
WxWVxxx\\
WxWWxxx\\
WxWxVxx\\
WxWxWxx
\end{array} \]
\end{minipage}\\ \\

Because of lack of space, in $T_3$ figure only the first two and
the last two columns, the four columns in the middle having being
omitted. After labeling these word expressions from $1$ to $40$,
the tableaux can be perused \footnote{The importance of perusal
and inspection of tables was aptly pointed out by D.H. Lehmer in
hia article MAA Studies in Mathematics, Vol 6, 1969} easily as
seen below.

\[ \begin{array}{c}
\underline{\text{$T_{1}$}}\\\\
1\\
2\\ \\ \\ \\ \\ \\
\end{array}
\,\,\,\,\,\,\,\,\,\, \begin{array}{c}
\underline{\text{$T_{2}$}}\\\\
1\,\,\,2\\
3\,\,\,4\\
5\,\,\,6\\
7\,\,\,8\\ \\ \\ \\
\end{array}\,\,\,\,\,\,\,\,\,\,
\begin{array}{ccccccccc}
 & & & & \underline{\text{$T_{3}$}}\\\\
1&2&3&4&5&6&7&8 \\
9&10&11&12&13&14&15&16 \\
17&18&19&20&21&22&23&24 \\
25&26&27&28&29&30&31&32 \\
5&6&13&14&33&34&35&36 \\
7&8&15&16&37&38&39&40 \\
\\
\end{array} \]

\par\vspace{10pt}

The general tableau of order $n$ has $2n$ lines and
$2^{n-1}\,S_{n-1}$ columns, that is a total of
$2^{n}\,n\,S_{n-1}$ entries. For $n\geq3$ it is easy to see that
some $n-$iterates appear in tableau $T_n$ with multiplicities
higher than $1$, as can be verified in tableau $T_3$. To prove it
we have to show that $2^n\,n\,S_{n-1}>2^n\,S_n$ for $n\geq3$. The
easy proof is as follows. Using the recursion
$S_n=\dfrac{2(2n-1)}{n+1}S_{n-1}$ for the Catalan numbers we have

\begin{align*}
2^n\,n\,S_{n-1}&>2^n\,S_n \\
n\,S_{n-1}&>\dfrac{2(2n-1)}{n+1}S_{n-1} \\
n(n+1)&>2(2n-1).
\end{align*}

The last inequality being true for $n\geq3$, application of the
pigeonhole principle does the rest.
\par\vspace{10pt}All concepts
and definitions of \cite{Petridi1} relating to one binary
operation $V(x,\,y)$ can be carried over literally to the present
case. Regrettably, because of the reasons explained in section
$1$, we were unable to go further than tableau $T_3$. We were
lucky, however, to discover that already for the incidence matrix
of this tableau the fundamental theorem of subsection $2.4$ of
\cite{Petridi1}, which is the key enabling to calculate the
number $I_n$ of formally reducible identities, remains true.
Because of the highly peculiar nature of this property we surmise
that it is equally true for all higher tableaux $T_n$. For easy
reference we repeat the theorem hereunder.
\par\vspace{10pt}
\begin{theorem}
Let \begin{enumerate}
      \item $\delta(J^{n}_{i},\,J^{n}_{j})=\left\{\begin{array}{c} \hspace{-0.2cm}1\,\text{if}\,\,J^{n}_{i}=J^{n}_{j}\,\,\text{reducible} \\
                                                                   0\,\text{if}\,\,J^{n}_{i}=J^{n}_{j}\,\,\text{irreducible}\end{array}\right.$
      \item $M(J_{i}^{n})=\,\,\text{the multiplicity of}\,\,J_{i}^{n}\,\,\text{in tableau}\,\,A_n$
      \item $I_n= \sum_{i,\,j}^{S_n}\delta(J^{n}_{i},\,J^{n}_{j})=\,\,\text{the number of
      reducible}\,\,n-\text{identities}$
       \item $\sum_{j=1}^{S_n}\delta(J^{n}_{i},\,J^{n}_{j})= \text{the number of reducible}
       \,\,n-\text{identities on the}\,\,i-\text{th}$ \\ $\text{line  of the incidence matrix of
       tableau}\,\,A_n,$
\end{enumerate}
then
$$\sum_{j=1}^{S_n}\delta(J^{n}_{i},\,J^{n}_{j})=\sum_{\nu=1}^{M(J_{i}^{n})}(-1)^{\nu-1}
\binom{M(J_{i}^{n})}{\nu}S_{n-\nu}.$$

\end{theorem}

Expressed in words the theorem says that the number of reducible
$n-$identities on the $i-$line of the incidence matrix of tableau
$A_n$ does not depend on $J^{n}_{i}$ but only on its multiplicity
$M(J_{i}^{n})$. An immediate consequence is that

\[ I_n=\sum_{k=1}^{[\frac{n+1}{2}]}T_{nk}\big(\sum_{\nu=1}^{k}(-1)^{\nu-1}\binom{k}{\nu}S_{n-\nu}\big)
\]

where $T_{nk}$ is the number of iterates in tableau $T_{n}$
having multiplicity $k$. As proved in \cite{Petridi1} $I_{n}$ is
$$I_{n}=\text{o}(1-e^{-\frac{n}{16}})$$
and the scarcity of the reducible identities is evinced by
$$S_{n}^{2}-I_{n}=\text{o}(e^{-\frac{n}{26}}).$$
\par\vspace{10pt}We now will prove the truth of this theorem for tableau $T_3$.
To this end we have calculated the incidence matrix relative to
tableau $T_3$ as shown in Exhibit attached hereto.
\par\vspace{10pt}
The proof leaps to the eye. Indeed, the four iterates
$J^{3}_{5}=5,\,J^{3}_{6}=6,\,J^{3}_{7}=7,\,J^{3}_{8}=8$ have all
multiplicity $2$, giving a sum $\sum\,1=14$. Similarly for the
iterates $J^{3}_{13},\,J^{3}_{14},\,J^{3}_{15},\,J^{3}_{16}$. All
other $32$ iterates have multiplicity $1$ with a sum $\sum\,1=8$.
Hence the number $I_{3}^{VW}$ of reducible $3-$identities is

\[ I_{3}^{VW}=32\cdot 8+8\cdot 14=368 \]

and the relative frequency is

\[ \dfrac{I_{3}^{VW}}{(2^3\,S_{3})^2}=\dfrac{368}{1600}=0.28. \]

\par\vspace{10pt}The main objective is of course to prove that

\[ \lim_{n\to\infty}\dfrac{I_{n}^{VW}}{(2^n\,S_{n})^2}=1, \]

which would imply that \textbf{irreducible} identities are
getting scarce with increasing $n$. Expressed otherwise, this
would mean that there are no algebras defined by \textbf{lengthy}
identities involving two binary operations.
\par\vspace{5pt}
The case of algebras with more than two binary operations
$V_1,\,V_2,\,\cdots,\,V_\lambda,\,\lambda>2$,
can be dealt with in the same way we did for $\lambda=2$.
The functional equation for the generating function
$\phi_{V_1,\,V_2,\,\cdots,\,V_\lambda}(t)$ turns out to be

\[ \dfrac{\phi_{V_1\,V_2\,\cdots\,V_\lambda}(t)-1}{t}=\lambda(\phi_{V_1\,V_2\,\cdots\,V_\lambda}(t))^{2} \]

which after solving gives the corresponding $``$Catalan" numbers
of the structure

\[ S_{n}^{V_1\,V_2\,\cdots\,V_\lambda}=\lambda^{n}S_{n}\,\,\;(S_{n}=n-\text{th Catalan number}). \]

\par\vspace{5pt}

\textbf{4 Operations of any arity}
\par\vspace{10pt}A direct approach to form the tableaux for several operations of arities
higher than two is well nigh impossible without powerful computer programs.
If at all even then. We may circumvent, however, the obstacle by reducing the problem to
the binary case as follows.
\par
\vspace{5pt} Given the operations
\begin{gather*}
V^{\alpha_1}(x_1,\,x_2,\,\cdots,\,x_{\alpha_1}) \\
V^{\alpha_2}(x_1,\,x_2,\,\cdots,\,x_{\alpha_2}) \\
\cdots\cdots\cdots\cdots\cdots\cdots \\
V^{\alpha_m}(x_1,\,x_2,\,\cdots,\,x_{\alpha_m})
\end{gather*}
we form their respective $\binom{\alpha_{i}}{2}$, projections on
the subspaces of the variable places
$$
V_{jkj}^{\alpha_1}(x,\,y)=V^{\alpha_1}(\underbrace{c,\,\cdots,\,c}_{i},\,x,\,\underbrace{c,\,\cdots,\,c}_{j}\,x\,\underbrace{c,\,\cdots,\,c}_{k})
,\,$$ taken over all solutions of  $i+j+k=\alpha_1-2,$
$$
V_{jkj}^{\alpha_2}(x,\,y)=V^{\alpha_2}(\underbrace{c,\,\cdots,\,c}_{i},\,x,\,\underbrace{c,\,\cdots,\,c}_{j}\,x\,\underbrace{c,\,\cdots,\,c}_{k})
,\,$$ taken over all solutions of  $i+j+k=\alpha_2-2,$
$$
\cdots\,\,\cdots\,\,\cdots\,\,\cdots\,\,\cdots\,\,\cdots\,\,\cdots\,\,\cdots\,\,\cdots\,\,\cdots\,\,\cdots\,\,\cdots\,\,\cdots\,\,\cdots\,\,\cdots\,\,\cdots
\,\,\cdots\,$$
$$
V_{jkj}^{\alpha_m}(x,\,y)=V^{\alpha_m}(\underbrace{c,\,\cdots,\,c}_{i},\,x,\,\underbrace{c,\,\cdots,\,c}_{j}\,x\,\underbrace{c,\,\cdots,\,c}_{k})
,\,$$

taken over all solutions of  $i+j+k=\alpha_m-2.$
\par\vspace{10pt}
Since all these projections are binary operations we can apply to
them the results of the previous section $3$ and conclude that the
fundamental theorem holds true.

\par\vspace{20pt}

\textbf{5 Connection with Formal Languages}

\par\vspace{20pt}Seen from the angle of Formal Languages, a set of
operations and their iterates is just the Language $L$ generated
by the grammar
$G(V_{1}^{\alpha_{1}},V_{2}^{\alpha_{2}},\ldots,V_{k}^{\alpha_{k}},x)$
with $$x\in L:\;\;\text{the starting word}$$ and the derivation
rules of words
\par
 $$\text{If}\hspace{1cm}W_{1}\in L$$ $$\text{If}\hspace{1cm}W_{2}\in L$$ $$\ldots$$
$$\text{If}\hspace{1cm}W_{k}\in L$$
$$\text{Then}\hspace{1cm}V_{1}^{\alpha_{1}}W_{x_{1}}W{x_{2}}
\ldots W_{x_{\alpha_{1}}} \in L$$
 $$\text{Then}\hspace{1cm}V_{2}^{\alpha_{2}}W_{y_{1}}W{y_{2}} \ldots
W_{y_{\alpha_{2}}} \in L$$ $$\ldots  \;\;  \ldots$$
 $$\text{Then}\hspace{1cm}V_{k}^{\alpha_{k}}W_{z_{1}}W{z_{2}} \ldots
W_{z_{\alpha_{k}}} \in L$$ where the indexes $x_{1},\ldots
x_{\alpha_{1}},\; y_{1},\ldots y_{\alpha_{2}},\ldots z_{1},\ldots
z_{\alpha_{k}}$, run over all permutations with repetitions of
$\{ 1,2,\ldots, k\}$.
\par
\vspace{5pt} The reverse is also true. If in the alphabet of the
grammar all non-terminal symbols are replaced by $x$ and the
terminal symbols are replaced respectively by
$V_{1}^{\alpha_{1}},V_{2}^{\alpha_{2}},\ldots V_{k}^{\alpha_{k}}$
we obtain the structure with operations
$V_{1}^{\alpha_{1}},V_{2}^{\alpha_{2}},\ldots V_{k}^{\alpha_{k}}$.
\par\vspace{5pt} For further reading on the subject see
sub-section 4.2 of \cite{Petridi1}. Whether precise analytical
results, analogous to those of \cite{Petridi1} and the present
paper, hold for all Formal Languages as well as their uses in
Information Theory is, to our knowledge, an open field to be
explored.
\newpage
\begin{Exhibit}[H]
\begin{center}
\textbf{Exhibit}\\
Incidence matrix $||\delta(J^{3}_{i},\,J^{3}_{j})||$ relative to $T_3$ \\
($J^{3}_{i}$ is denoted by $i$. Blanc spaces mean $0'$s)\\
\end{center}
\[\hspace{-1.5cm}\footnotesize \arraycolsep=2pt
{\begin{array}{ccccccccccccccccccccccccccccccccccccccccccccc}
i\backslash j&&1&2&3&4&5&6&7&8&9&10&1&2&3&4&5&6&7&8&9&20&1&2&3&4&5&6&7&8&9&30&1&2&3&4&5&6&7&8&9&40&\,\,\,\,\sum_{i}\,1&\,\,\,M(J^{3}_{i})\\
\\
1& &1&1&1&1&1&1&1&1& & & & & & & & & & & & & & & & & & & & & & & & & & & & & & & & &8&1\\
2& &1&1&1&1&1&1&1&1& & & & & & & & & & & & & & & & & & & & & & & & & & & & & & & & &8&1\\
3& &1&1&1&1&1&1&1&1& & & & & & & & & & & & & & & & & & & & & & & & & & & & & & & & &8&1\\
4& &1&1&1&1&1&1&1&1& & & & & & & & & & & & & & & & & & & & & & & & & & & & & & & & &8&1\\
5& &1&1&1&1&1&1&1&1& & & & & 1&1& & & & & & & & & & & & & & & & & & &1&1&1&1& & & & &14&2\\
6& &1&1&1&1&1&1&1&1& & & & & 1&1& & & & & & & & & & & & & & & & & & &1&1&1&1& & & & &14&2\\
7& &1&1&1&1&1&1&1&1& & & & & & &1&1& & & & & & & & & & & & & & & & & & & & &1&1&1&1&14&2\\
8& &1&1&1&1&1&1&1&1& & & & & & &1&1& & & & & & & & & & & & & & & & & & & & &1&1&1&1&14&2\\
9& & & & & & & & & & 1&1&1&1&1&1&1&1& & & & & & & & & & & & & & & & & & & & & & & & &8&1\\
10& & & & & & & & & & 1&1&1&1&1&1&1&1& & & & & & & & & & & & & & & & & & & & & & & & &8&1\\
11& & & & & & & & & & 1&1&1&1&1&1&1&1& & & & & & & & & & & & & & & & & & & & & & & & &8&1\\
12& & & & & & & & & & 1&1&1&1&1&1&1&1& & & & & & & & & & & & & & & & & & & & & & & & &8&1\\
13& & & & & &1&1& & & 1&1&1&1&1&1&1&1& & & & & & & & & & & & & & & & &1&1&1&1& & & & &14&2\\
14& & & & & &1&1& & & 1&1&1&1&1&1&1&1& & & & & & & & & & & & & & & & &1&1&1&1& & & & &14&2\\
15& & & & & & & &1&1& 1&1&1&1&1&1&1&1& & & & & & & & & & & & & & & & & & & & &1&1&1&1&14&2\\
16& & & & & & & &1&1& 1&1&1&1&1&1&1&1& & & & & & & & & & & & & & & & & & & & &1&1&1&1&14&2\\
17& & & & & & & & & & & & & & & & & & 1&1&1&1&1&1&1&1& & & & & & & & & & & & & & & & &8&1\\
18& & & & & & & & & & & & & & & & & & 1&1&1&1&1&1&1&1& & & & & & & & & & & & & & & & &8&1\\
19& & & & & & & & & & & & & & & & & & 1&1&1&1&1&1&1&1& & & & & & & & & & & & & & & & &8&1\\
20& & & & & & & & & & & & & & & & & & 1&1&1&1&1&1&1&1& & & & & &
& & & & & & & & & & &8&1\\
21& & & & & & & & & & & & & & & & & & 1&1&1&1&1&1&1&1& & & & & & & & & & & & & & & & &8&1\\
22& & & & & & & & & & & & & & & & & & 1&1&1&1&1&1&1&1& & & & & & & & & & & & & & & & &8&1\\
23& & & & & & & & & & & & & & & & & & 1&1&1&1&1&1&1&1& & & & & & & & & & & & & & & & &8&1\\
24& & & & & & & & & & & & & & & & & & 1&1&1&1&1&1&1&1& & & & & & & & & & & & & & & & &8&1\\
25& & & & & & & & & & & & & & & & & & & & & & & & & & 1&1&1&1&1&1&1&1& & & & & & & & &8&1\\
26& & & & & & & & & & & & & & & & & & & & & & & & & & 1&1&1&1&1&1&1&1& & & & & & & & &8&1\\
27& & & & & & & & & & & & & & & & & & & & & & & & & & 1&1&1&1&1&1&1&1& & & & & & & & &8&1\\
28& & & & & & & & & & & & & & & & & & & & & & & & & & 1&1&1&1&1&1&1&1& & & & & & & & &8&1\\
29& & & & & & & & & & & & & & & & & & & & & & & & & & 1&1&1&1&1&1&1&1& & & & & & & & &8&1\\
30& & & & & & & & & & & & & & & & & & & & & & & & & & 1&1&1&1&1&1&1&1& & & & & & & & &8&1\\
31& & & & & & & & & & & & & & & & & & & & & & & & & & 1&1&1&1&1&1&1&1& & & & & & & & &8&1\\
32& & & & & & & & & & & & & & & & & & & & & & & & & & 1&1&1&1&1&1&1&1& & & & & & & & &8&1\\
33& & &1&1& & & & & & &1&1& & & & & & & & & & & & & & & & & & & & & &1&1&1&1& & & & &8&1\\
34& & &1&1& & & & & & &1&1& & & & & & & & & & & & & & & & & & & & & &1&1&1&1& & & & &8&1\\
35& & &1&1& & & & & & &1&1& & & & & & & & & & & & & & & & & & & & & &1&1&1&1& & & & &8&1\\
36& & &1&1& & & & & & &1&1& & & & & & & & & & & & & & & & & & & & & &1&1&1&1& & & & &8&1\\
37& & & & &1&1&& & & & & &1&1&& & & & & & & & & & & & & & & & & & & & & & &1&1&1&1&8&1\\
38& & & & &1&1&& & & & & &1&1&& & & & & & & & & & & & & & & & & & & & & & &1&1&1&1&8&1\\
39& & & & &1&1&& & & & & &1&1&& & & & & & & & & & & & & & & & & & & & & & &1&1&1&1&8&1\\
40& & & & &1&1&& & & & & &1&1&& & & & & & & & & & & & & & & & & & & & & & &1&1&1&1&8&1

\end{array}}\]
\end{Exhibit}
\par
From above table we get
$$I_{3}^{VW}=\sum_{i=1}^{40}\sum_{j=1}^{40}1 = 368$$

\end{document}